\input amssym.def
\input epsf

\let \blskip = \baselineskip
\parskip=1.2ex plus .2ex minus .1ex

\tabskip 20pt
\tolerance = 1000
\pretolerance = 50
\newcount\itemnum
\itemnum = 0
\overfullrule = 0pt

\def\title#1{\bigskip\centerline{\bigbigbf#1}}
\def\author#1{\bigskip\centerline{\bf #1}\smallskip}
\def\address#1{\centerline{\it#1}}
\def\abstract#1{\vskip1truecm{\narrower\noindent{\bf Abstract.} #1\bigskip}}

\def\sp{\bigskip}
\def\nosp{\vskip -\the\blskip plus 1pt minus 1pt}

\def\br{\hfil\break} 
\def\ti{\br \hglue \the \parindent}

\def\ce#1{\LP\centerline{#1}}

\def\skipit#1{}
\def\mdag{\raise 3pt\hbox{\dag}}

\def\XP{\par\noindent\hang}
\def\LP{\par\noindent}
\def\BP[#1]{\par\item{[#1]}}
\def\SH#1{\sp\vskip\parskip\leftline{\bigbf #1}\nobreak}

\def\TH#1{\sp\XP{\bf THEOREM\ \shead#1}}
\def\LM#1{\sp\XP{\bf LEMMA\ \shead#1}}

\def\CO#1{\sp\XP{\bf COROLLARY\ \shead#1}}

\def\RK#1{\sp\XP{\bf REMARK\ \shead#1}}

\def\PF{\LP{\bf Proof:\ }}
\def\NX{\advance\itemnum by 1 \sp\LP {\bf \shead \the\itemnum.\ }}
\def\qed{\null\nobreak\hfill\hbox{${\vrule width 5pt height 6pt}$}\par\sp}

\def\cart{\>\hbox{${\vcenter{\vbox{
    \hrule height 0.4pt\hbox{\vrule width 0.4pt height 4.5pt
    \kern4pt\vrule width 0.4pt}\hrule height 0.4pt}}}$}\>}
\def\bxmu{\>\hbox{${\vcenter{\vbox {
    \hrule height 0.4pt\hbox{\vrule width 0.4pt height 4pt
    \hskip -1.3pt\lower 1.8pt\hbox{$\times$}\negthinspace\vrule width 0.4pt}
    \hrule height 0.4pt}}}$}\>}

\def\lin#1{\hbox to #1true in{\hrulefill}}


\def\foot#1{\raise 6pt \hbox{#1} \kern -3pt}

\def\fig #1 #2 #3 #4 #5 {\sp \ce{ {\epsfbox[#1 #2 #3 #4]{figs/#5.ps}} }}



\def\JCTB{{\it J.\ Comb.\ Theory (B)}}

\def\JGT{{\it J.\ Graph Theory}}

\def\DM{{\it Discrete Math.{}}}

\def\al{\alpha}			
    
\def\sg{\sigma}	 \def\GA{\Gamma}



\def\nul{\hbox{\O}}    
		
\def\esub{\subseteq}

	\def\join{\lor}	


\def\({\left(}	\def\){\right)}

\def\small#1{{\scriptstyle #1}}

\def\CH#1#2{{{#1}\choose{#2}}}

\def\tFR#1#2{\small{{#1 \over #2}}}

\def\FL#1{\left\lfloor{#1}\right\rfloor}

\def\UM#1#2{\bigcup_{#1\in#2}}

\def\VEC#1#2#3{#1_{#2},\ldots,#1_{#3}}
\def\VECOP#1#2#3#4{#1_{#2}#4 \cdots #4 #1_{#3}}

\def\st{\colon\;} 

\def\SET#1:#2{\{#1\colon\;#2\}}

\def\B#1{{\bf #1}}		
\def\C#1{\left | #1 \right |}    

\def\bA{\B A}  \def\bC{\B C}

  \def\bR{\B R}

\def\Pb{\ov{P}}



  \def\Koe{K_{1,3}} 




\magnification=\magstep1
\vsize=9.0 true in
\hsize=6.5 true in
\headline={\hfil\ifnum\pageno=1\else\folio\fi\hfil}
\footline={\hfil\ifnum\pageno=1\folio\else\fi\hfil}

\parindent=20pt
\baselineskip=12pt
\parskip=.5ex  

\def\shead{ }

\font\bigbf = cmb10 scaled \magstep1

\font\bigbigbf = cmb10 scaled \magstep2


\def\gpic#1{#1
     \bigskip\par\noindent{\centerline{\box\graph}}
     \medskip} 

\def\title#1{\bigskip\centerline{\bigbigbf#1}}
\def\author#1{\bigskip\centerline{\bf #1}}
\def\address#1{\centerline{\it#1}}
\def\abstract#1{\vskip10pt{\narrower\noindent{\bf Abstract.} #1\bigskip}}


\title{THE LEAFAGE OF A CHORDAL GRAPH}
\author{In-Jen Lin} 
\address{National Ocean University, Taipei, Taiwan}
\author{Terry A. McKee\foot{1}}
\address{Wright State University, Dayton, OH 45435-0001}
\author{Douglas B. West\foot{2}}
\address{University of Illinois, Urbana, IL 61801-2975}
\vfootnote{}{\br
   \foot{1}Research supported in part by ONR Grant N00014-91-J-1210.\br
   \foot{2}Research supported in part by NSA/MSP Grant MDA904-93-H-3040.\br
   Running head: LEAFAGE OF CHORDAL GRAPHS\br
   Keywords: chordal graph, subtree intersection representation, leafage\br
   Written July 1994, revised December 1996 and April 1998.}
\abstract{
  The {\it leafage} $l(G)$ of a chordal graph $G$ is the minimum number of
  leaves of a tree in which $G$ has an intersection representation by subtrees.
  We obtain upper and lower bounds on $l(G)$ and compute it on special classes.
  The maximum of $l(G)$ on $n$-vertex graphs is $n-\lg n-\tFR12\lg\lg n+O(1)$.
  The {\it proper leafage} $l^*(G)$ is the minimum number of leaves when no
  subtree may contain another; we obtain upper and lower bounds on $l^*(G)$.
  Leafage equals proper leafage on claw-free chordal graphs.
  We use asteroidal sets and structural properties of chordal graphs.}

\def\ast{asteroidal}
\def\peo{perfect elimination ordering}
\def\rep{representation}
\def\ind{independent}
\def\intn{intersection}
\def\intg{intersecting}
\def\half{\tFR12}
\def\charz{characterization}

\def\Pb{\bar P}

\SH
{1. INTRODUCTION}
A simple graph is {\it chordal} (or {\it triangulated}) if every
cycle of length exceeding 3 has a chord.  The {\it \intn\ graph} of a family of
sets is the graph defined by assigning one vertex for each set and joining two
vertices by an edge if and only if the corresponding sets intersect.  A graph is
chordal if and only if it is the \intn\ graph of a family of subtrees of a host
tree [3,8,23]; such a family is a {\it subtree \rep} of the graph.  The {\it
interval graphs} are the chordal graphs having subtree \rep s in which the host
tree is a path; this allows the subtrees to be viewed as intervals on the real
line.  Given the many applications of interval graphs and the ease of
computation on interval graphs, it is natural to ask for measures of how far a
chordal graph is from being an interval graph.

The {\it leafage} $l(G)$ of a chordal graph $G$ is the minimum number of leaves
of the host tree in a subtree \rep\ of $G$.  Interval graphs are the chordal
graphs with leafage at most two.  We derive bounds on leafage and study classes
in which equality holds, including $k$-trees, block graphs, and chordal graphs
having a dominating clique.  We prove that the maximum leafage of a chordal
graph on $n$ vertices is the maximum $k$ such that
$k\le\CH{n-k}{\FL{(n-k)/2}}$; this is $n-\lg n-\half\lg\lg n+O(1)$.  Our proofs
provide algorithms for computing leafage in some special classes.

Among the interval graphs are the {\it proper interval graphs}, which are those
representable using a family of intervals in which no interval properly contains
another.  As observed in [8], every chordal graph has a subtree \rep\ in which
no subtree properly contains another; call this a {\it proper subtree \rep}.  By
analogy with leafage, the {\it proper leafage} $l^*(G)$ of a chordal graph $G$
is the minimum number of leaves in a proper subtree \rep of $G$.  We extend a
\charz\ of proper interval graphs by Roberts [19] to obtain a lower bound on
proper leafage: $l^*(G)$ is at least the maximum number of ``modified extreme
points'' occurring in an induced subgraph of $G$.  We present examples where
the bound is arbitrarily bad, but the bound is sharp for $k$-trees, where it
equals the number of simplicial vertices.

A notion analogous to leafage has already been characterized for a \rep\
model based on containment.  Leclerc [14] proved that a graph is the
comparability graph of a partial order of dimension at most $d$ if and only if
it is the containment poset of a family of subtrees of a tree with at most $d$
leaves.  Our notion of leafage has been applied in [17], and the analogue of
leafage for directed graphs is studied in [16].

We use well-known properties of chordal graphs from [3,6,11,20,23].  A {\it
simplicial vertex} in a graph is a vertex $v$ whose neighborhood $N(v)$ induces
a clique.  A {\it cutset} (also called {\it separating set} or
{\it vertex cut}) is a set of vertices whose deletion leaves a disconnected
subgraph.  Every cutset of a chordal graph induces a clique, and thus by
induction every non-clique chordal graph has a nonadjacent pair of simplicial
vertices [6].  The neighborhood of a simplicial vertex is a {\it simplicial
neighborhood}.  The subgraph of $G$ induced by the set $S(G)$ of simplicial
vertices is a disjoint union of cliques.  We use $S(G)$ to denote this subgraph
in the same way that $Q$ may denote a clique or its set of vertices.  The
deletion of a simplicial vertex cannot affect the existence of chordless cycles.
Hence another \charz\ of a chordal graph [7,20] is the existence of a {\it
perfect elimination ordering}, meaning a vertex ordering $\VEC vn1$ in which the
vertices can be deleted such that for each $i$, $N(v_i)\cap\{\VEC vi1\}$ induces
a clique.

Given a chordal graph $G$, the {\it derived graph} $G'$ is the induced subgraph
obtained by deleting all the simplicial vertices of $G$.  As an induced subgraph
of a chordal graph, $G'$ also is a chordal graph.  We describe a subtree of a
tree by the set of vertices inducing it, and with this understanding we write
$f(v)$ for the subtree representing $v$ in a subtree \rep\ $f$.  We use $m(T)$
for the number of leaves of a tree $T$.

Every pairwise \intg\ family of subtrees of a tree has a common vertex [11,
p92].  Hence if $Q$ is a clique in $G$, any subtree \rep\ $f$ of $G$ assigns
some host vertex to all of $Q$.  If distinct vertices $q,q'$ are assigned to
cliques $Q,Q'$ with $Q\subset Q'$, then for $v\in Q$ the entire $q,q'$-path in
the host belongs to $f(v)$.  The first edge on this path can be contracted to
obtain a smaller \rep\ without changing the number of leaves.  Therefore, we
may restrict our attention to {\it minimal \rep s}, which are those subtree
\rep s having a bijection between the maximal cliques of $G$ and the vertices
of the host tree.  This restriction is {\it not valid} for proper leafage.  We
use the term {\it optimal representation} to refer to a subtree representation
having the minimum number of host leaves (subject to appropriate conditions).

Gavril [9] and Shibata [22] proved that minimal representations correspond
to maximum weight spanning trees in the weighted clique graph of $G$, where the
{\it weighted clique graph} has a vertex for each maximal clique of $G$, and the
weight of the edge $QQ'$ is $\C{Q\cap Q'}$.  Thus $l(G)$ equals the minimum
number of leaves in a maximum weight spanning tree of the weighted clique graph
of $G$.  This observation does not seem to simplify the problem.

\SH
{2. ASTEROIDAL SETS AND SPECIAL CLASSES}
Our first lower bound for $l(G)$ generalizes the notion of asteroidal triple.
We prove constructively that this bound is exact for trees and more generally
for $k$-trees.

An {\it asteroidal triple} in a graph is a triple of distinct vertices such that
each pair is connected by some path avoiding the neighborhood of the third
vertex.  Interval graphs are the chordal graphs without asteroidal triples
([10,15]).  A set $S\esub V(G)$ is an {\it asteroidal set} if every triple of
vertices in $S$ is an asteroidal triple.  This concept arises also in [18,23].
Let $a(G)$ denote the maximum size of an \ast\ set in $G$.  This parameter
has subsequently been named the {\it asteroidal number} in [1,2,12,13]; these
papers explore further properties of asteroidal sets, asteroidal number,
and leafage.

We define a notion like asteroidal sets for subtrees.  If $T_i$, $T_j$, $T_k$
are subtrees of a tree, then $T_k$ is {\it between} $T_i$ and $T_j$ if $T_i$ and
$T_j$ are disjoint and the unique path connecting them intersects $T_k$.  A
collection of pairwise disjoint subtrees such that none is between two others is
an {\it asteroidal collection} of subtrees.

\LM 1.
If $\VEC T1n$ is an asteroidal collection of subtrees of a tree $T$,
then $T$ has at least $n$ leaves.
\PF
For each leaf $v$ of $T$, we assign to $v$ the first subtree in $\{T_i\}$
encountered on the path from $v$ to the nearest vertex of degree at least 3.
Such a subtree exists, else we could delete the edges of that path to reduce
$m(T)$ without affecting $\{T_i\}$.  If $m(T)<n$, then some subtree $T_k$ in our
list is assigned to no leaf.  Let $x$ be a vertex of $T_k$, and let $P$ be a
maximal path in $T$ containing $x$.  The endpoints of $P$ are leaves of the
tree, and $T_k$ is between the subtrees assigned to those leaves.  The
contradiction yields $m(T)\ge n$.  \qed

\TH 1.
In a subtree \rep\ of a chordal graph $G$, the subtrees corresponding to an
asteroidal set of vertices form an asteroidal collection, and $l(G)\ge a(G)$.
\PF
Let $f$ be a subtree \rep\ of $G$ in $T$.  Every asteroidal triple of vertices
is an \ind\ set, and thus every asteroidal set is an \ind\ set.  Hence the
subtrees corresponding to an \ast\ set in a subtree \rep\ of $G$ are pairwise
disjoint.  If none is between two others, then the subtrees form an asteroidal
collection in $T$, and Lemma 1 yields the bound.

Suppose that $f(w)$ is between $f(u)$ and $f(v)$.  Let $P$ be a $u,v$-path in
$G$ containing no neighbor of $w$.  Because every two successive vertices in $P$
are adjacent, the corresponding subtrees in $T$ have a common point, and hence
the union $\UM xP f(x)$ of the subtrees representing these vertices is
a connected subgraph $T^-$ of $T$.  Since $P$ contains no neighbor of $w$, $T^-$
must be disjoint from $f(w)$, but this contradicts the assumption that the
unique path between $f(u)$ and $f(v)$ in $T$ contains a vertex of $f(w)$.  \qed

\CO 1.
If $G$ is a tree other than a star, then $l(G)$ is the number of leaves in the
derived tree $G'$.
\PF
Suppose that $G'$ has $k$ leaves.  For the lower bound, we obtain a asteroidal
set.  For each leaf of $G'$, select a leaf of $G$ adjacent to it.  This yields
an asteroidal set $S$ of size $k$, because for $x,y \in S$, the unique
$x,y$-path in $G$ consists of $x$, $y$, and a path between leaves of $G'$ that
are adjacent to $x$ and $y$.  This path contains no other leaf of $G$ or $G'$,
so it avoids the neighborhoods of other vertices in $S$.

For the upper bound, we construct a subtree \rep.  For the host tree $T$,
begin with the tree obtained from $G'$ by subdividing each edge once.  Next, for
each $v'\in V(G')$, choose one edge incident to $v'$ in the current host and
subdivide it $k$ times, where $k$ is the number of leaves of $G$ incident
to $v'$.  Use each of these new vertices to represent one leaf of $G$ incident
to $v'$, and let $f(v')$ consist of the vertex $v'$ in $T$ together with the
vertices introduced by subdividing edges incident to $v'$.  This yields a
\rep\ of $G$.  Since $G$ is not a star, $G'$ has at least two leaves, and $T$
and $G'$ have the same number of leaves.  \qed

A $k$-{\it tree} is a chordal graph that can be constructed from $K_k$ by a
sequence of vertex additions in which the neighborhood of each new vertex is a
$k$-clique of the current graph.  With additional lemmas, we can generalize
the construction of Corollary 1 to $k$-trees.

\LM 2.
In a non-clique $k$-tree, the simplicial vertices are the vertices of degree
$k$ and form an independent set.
\PF
The smallest non-clique $k$-tree has $k+2$ vertices and satisfies the claim.
When a vertex receives a new neighbor, its degree exceeds $k$, and its old
neighborhood cannot be contained among the other $k-1$ neighbors of the new
simplicial vertex.  Hence vertices of degree exceeding $k$ (and neighbors
of simplicial vertices) are not simplicial.  \qed

Every minimal cutset of a $k$-tree induces a $k$-clique (Rose [20]).  Indeed,
non-clique $k$-trees are precisely the connected graphs in which the largest
clique has $k+1$ vertices and every minimal cutset induces a $k$-clique.  Since
deleting a simplicial vertex of a $k$-tree leaves a smaller $k$-tree, the
construction procedure defining a $k$-tree can begin from any $k$-clique.
These properties yields the following statement, which does not hold
for general chordal graphs (it fails for the interval graph $2P_4\join K_1$).
(The {\it join} $G\join H$ of two graphs $G$ and $H$ is obtained from the
disjoint union $G+H$ by adding as edges $\{uv\st u\in V(G), v\in V(H)\}$.)

\LM 3.
If $G$ is a non-clique $k$-tree having distinct simplicial neighborhoods, then
$G$ has distinct simplicial neighborhoods that are not cutsets of $G'$.
\PF
By Lemma 2, every simplicial neighborhood in $G$ is contained in $V(G')$.
We call a simplicial vertex $v$ {\it good} in $G$ if $G'-N(v)$ is connected.

If $G$ is a non-clique $k$-tree, then $G$ has only one simplicial neighborhood
if and only if $G'=K_k$.  If $G'\ne K_k$, then $G$ has a simplicial vertex
$x$ such that $G-x$ has distinct simplicial neighborhoods unless $G$ has $k+3$
vertices and $G'=K_{k+1}$.  There is one such $k$-tree, and its two simplicial
vertices are good.  This serves as the basis for induction.

For the induction step, choose $x\in S(G)$.  The induction hypothesis implies
that $G-x$ has two good simplicial vertices $u,v$ with distinct neighborhoods.
If $G'=(G-x)'$, then every simplicial vertex of $G-x$ is also simplicial in $G$.
In this case, $u,v$ are good for $G$.

Hence we may assume that $G'\ne(G-x)'$, which means that some simplicial vertex
of $G-x$ is not simplicial in $G$.  Such a vertex $y$ must be a neighbor of $x$.
Applied to $G-x$, Lemma 2 implies both that $y$ is unique and that 
$(G-x)'=G'-y$.

If $y\notin N(v)$, then $(G-x)'-N(v)=[G'-N(v)]-y$.  Hence $v$ remains good for
$G$ unless $y$ is isolated in $G'-N(v)$ (similarly for $u$).  Since $y$ and $v$
each have degree $k$ in $G-x$, this requires $N_{G-x}(y)=N_{G-x}(v)$, so we may
assume that $y=v$ and that $u$ is good in $G$.  Also, $u,x$ have distinct
neighborhoods in $G$, since $y\in N(x)$ and $y=v$ is simplicial in $G-x$.
Thus it suffices to prove that $x$ is good in $G$.

Let $z$ be the unique member of $N(y)-N(x)$.  Since $y$ is good in $G-x$ and
$(G-x)'-N(y) = [G'-N(x)]-z$, it suffices to show that $z$ is not isolated in
$G'-N(x)$.  If so, then $z$ has exactly $k$ neighbors in $G'$.  By Lemma 2,
$z\in S(G')$, but $G'$ cannot have adjacent simplicial vertices.  \qed

\skipit{
Let $Q$ be a $k$-clique separating $G$; as observed above,
$Q\cap S(G)=\nul$.  Every isolated vertex of $G-Q$ is simplicial in $G$.
If $G-Q$ is an independent set, then $G = (n-k)K_1\join K_k$ for some
$n\ge k+2$.  In this case every vertex of $G-Q$ is good, since $G'$ is a clique
and has no cutset.
   Lemma 2 guarantees that $G'$ contains the full neighborhood of any simplicial
vertex of $G$.  Say that $v$ is {\it good} if $v\in S(G)$ and $G'-N(v)$ is
connected; we use induction on the order of $G$ to guarantee nonadjacent good
vertices.  Let $Q$ be a $k$-clique separating $G$; as observed above,
$Q\cap S(G)=\nul$.  Every isolated vertex of $G-Q$ is simplicial in $G$.
If $G-Q$ is an independent set, then $G = (n-k)K_1\join K_k$ for some
$n\ge k+2$.  In this case every vertex of $G-Q$ is good, since $G'$ is a clique
and has no cutset.
   For the remaining cases, we will prove that if $F$ is a nontrivial component
of $G-Q$, then $F$ contains a good vertex of $G$.  If $G-Q$ has only one
nontrivial component $F$, we can then choose a good vertex from $F$ and choose
any isolated vertex $v$ of $G-Q$; deleting $N(v)=Q$ from $G'$ leaves a connected
graph that is a subgraph of $F$.  Otherwise $G-Q$ has at least two nontrivial
components, and we may choose a good vertex from any two of these.
   Let $F$ be a nontrivial component of $G-Q$, and let $H$ be the subgraph of
$G$ induced by $V(F)\cup Q$.  Since $F$ is nontrivial, $H$ is not a clique.
Since $V(F)$ has no edges to $G-V(H)$, starting from $Q$ as a root and
restricting the construction sequence to vertices of $F$ will grow $H$ from $Q$;
hence $H$ is a $k$-tree.  Since $Q$ is a clique, the induction hypothesis
guarantees the existence outside $Q$ of a good vertex $v$ for $H$.  Every
$v\in S(H)-Q$ is simplicial in both $H$ and $G$, since $V(F)$ has no edges to
$G-V(H)$.  In particular, $N_H(v)=N_G(v)\esub G'$.
   It suffices to show that $G'-N(v)$ is connected.  We have observed that the
subgraph of $G$ induced by $Q$ and the vertices of any component of $G-Q$ is a
$k$-tree that can be grown from $Q$; hence for each $x$ in $Q$ and $y$ outside
$Q$ in such a component, there is a path from $y$ to $Q$ that first reaches $Q$
at $x$. Since $H'-N(v)$ is connected and $Q\esub V(G')$, it thus suffices to
prove that $G'-N(v)$ contains a vertex of $Q$.  This is immediate, since
$Q\esub N(v)$ implies that $v$ is an isolated vertex of $G-Q$.  \qed
}

\TH 2.
If $G$ is a non-clique $k$-tree, then $l(G)=\max\{2,r(G)\}$, where $r(G)$ is the
number of distinct simplicial neighborhoods of $G$ that are not cutsets in $G'$.
\PF
When $G$ is a non-clique $k$-tree, $G$ has distinct simplicial neighborhoods
unless $G'$ is a nonempty clique, which occurs if and only if $G$ is the join of
a $k$-clique and an independent set.  Such a $k$-tree is an interval graph.

When $G$ has distinct simplicial neighborhoods, Lemma 3 implies that $r(G)\ge2$.
Let $\bR(G)$ be the set of simplicial neighborhoods of $G$ that are not cutsets
of $G'$.  For the lower bound, we construct an asteroidal set of size $r(G)$ by
selecting one simplicial neighbor of each $R\in\bR(G)$.  Given three such
vertices $x,y,z$, there is a $y,z$-path in $G$ avoiding $N(x)$ because
$G'-N(x)$ is connected.

For the upper bound, we use induction to construct a subtree \rep\ with the
desired number of leaves, in which for each $R\in\bR(G)$, the maximal cliques
containing the simplicial verties $\{v\in S(G)\st N(v)=R\}$ occur at a pendant
path in the host tree $T$, in any specified order.  We begin with the case
where $G'=K_k$.  Here $G$ is the join of $K_k$ and $j$ simplicial vertices,
and $T=P_j$, with each vertex of $G'$ assigned all of $T$ and each vertex of
$S(G)$ assigned one vertex of $T$.

For the induction step, we may assume that $G'$ has more than one maximal clique
and that $r(G)\ge2$.  By Lemma 3, we may choose $x\in S(G)$ such that $G'-N(x)$
is connected.  The graph $G-x$ is a smaller $k$-tree, and we apply the induction
hypothesis.  Let $f$ be a subtree \rep\ of $G-x$ in a host tree $T^-$ that has
the specified properties.  If $N(x)=N(y)$ for some $y\in S(G)$, then
$r(G-x)=r(G)$; here we obtain $T$ by inserting another vertex in the pendant
path in $T^-$ that corresponds to $N(x)$, assigning the new vertex to
$N(x)\cup\{x\}$.  If $N(x)\ne N(y)$ for all $y\in S(G)-\{x\}$, and no neighbor
of $x$ is simplicial in $G-x$, then $\bR(G-x)=\bR(G)-\{N(x)\}$, and we are
permitted to add a leaf in expanding $f$ to a subtree \rep\ of $G$.  Since
$N(x)$ is a clique in $G-x$, its corresponding subtrees have a common vertex in
$T^-$, and we can append a new leaf to that vertex, assigned to $N(x)\cup\{x\}$.

In the remaining case, some neighbor $z$ of $x$ is simplicial in $G-x$.  By
Lemma 2, $z$ is unique.  Let $Q=N_{G-x}(z)$, and let $q$ be the vertex assigned
by $f$ to the maximal clique $Q\cup\{z\}$.  If $Q$ is not a cutset of $(G-x)'$,
then by the induction hypothesis we may assume that $q$ is a leaf of $T^-$.  In
this case $Q$ is a cutset of $G'$ (it isolates $z$), so $r(G)=r(G-x)$ even if
$Q$ is a simplicial neighborhood in $G$.  On the other hand, if $Q$ is a cutset
of $(G-x)'$, then $r(G)=r(G-x)+1$.  In either case, we expand $T^-$ by adding a
leaf adjacent to $q$ and assign it to $N(x)\cup\{x\}$.  The resulting $T$ has
$r(G)$ leaves.  Furthermore, since $N(x)$ is not the neighborhood of any
simplicial vertex other than $x$, the claim about arbitrarily ordering the
cliques involving a given simplicial neighborhood also holds.  \qed

We prove that $l(G)= a(G)$ for one more class of chordal graphs.  The {\it
blocks} of a graph are its maximal subgraphs that have no cut-vertex.  A {\it
block graph} is a graph in which every block is a clique (equivalently, a
graph is a block graph if and only if it is the intersection graph of the blocks
in some graph).  Block graphs have no chordless cycles.  A {\it leaf block} of
$G$ is a block containing only one cut-vertex of $G$.

\TH 3.
If $G$ is a block graph that is not a clique, then $l(G) = \max\{2,r'(G)\}$,
where $r'(G)$ is the number of cut vertices of $G$ that are simplicial vertices
in $G'$.
\PF
In $G$, each maximal clique is a block, and a vertex belongs to $G'$ if and only
if it is a cut-vertex of $G$.  Also, if a cut vertex of $G$ is a simplicial
vertex of $G'$, then it belongs to a leaf block of $G$.  We form an asteroidal
set by selecting, for each simplicial vertex $v$ of $G'$, one simplicial vertex
from one leaf block of $G$ containing $v$.  This set has size $r'(G)$.  When
$r'(G)=1$, $G$ consists of two cliques sharing one vertex and is an interval
graph.

For the upper bound, we build by induction on $r'(G)$ a subtree \rep\ such that
the cliques that are leaf blocks containing a particular simplicial vertex of
$G'$ appear on a pendant path of the host tree.  If $r'(G)\le 2$, then $G'$ is a
path $\VEC v1p$.  To form the subtree \rep, first form a path $\VEC u0{2p+2}$
with $f(v_i) = \{u_{2i-2},u_{2i-1},u_{2i},u_{2i+1}\}$.  This assigns an edge
to each complete subgraph of $G'$.  Now subdivide edges as needed to insert a
vertex for each component of $S(G)$ in the path assigned to its neighborhood in
$G'$.

When $r'(G)>2$, we choose a simplicial vertex $v$ of $G'$.  Letting $U$ be the
set of neighbors of $v$ in leaf blocks of $G$, we have $r'(G-U)<r'(G)$.
The induction hypothesis provides a \rep\ for $G-U$; to any vertex of the host
assigned to clique containing $v$, we append a path having one vertex for each
leaf block of $G$ containing $v$.  These cliques appear at those vertices, and
the entire path is added to the subtree representing $v$.  \qed

Consider the chordal graph $G$ in Fig.~1 formed by adding one simplicial vertex
adjacent to each edge of $G'=2K_2\join K_1$.  Here $a(G)=3$ but $l(G)=4$, as we
will see.

\nosp
\gpic{
\expandafter\ifx\csname graph\endcsname\relax \csname newbox\endcsname\graph\fi
\expandafter\ifx\csname graphtemp\endcsname\relax \csname newdimen\endcsname\graphtemp\fi
\setbox\graph=\vtop{\vskip 0pt\hbox{%
    \graphtemp=.5ex\advance\graphtemp by 0.111in
    \rlap{\kern 0.778in\lower\graphtemp\hbox to 0pt{\hss $\bullet$\hss}}%
    \graphtemp=.5ex\advance\graphtemp by 0.496in
    \rlap{\kern 1.000in\lower\graphtemp\hbox to 0pt{\hss $\bullet$\hss}}%
    \graphtemp=.5ex\advance\graphtemp by 0.881in
    \rlap{\kern 0.778in\lower\graphtemp\hbox to 0pt{\hss $\bullet$\hss}}%
    \graphtemp=.5ex\advance\graphtemp by 0.881in
    \rlap{\kern 0.333in\lower\graphtemp\hbox to 0pt{\hss $\bullet$\hss}}%
    \graphtemp=.5ex\advance\graphtemp by 0.496in
    \rlap{\kern 0.111in\lower\graphtemp\hbox to 0pt{\hss $\bullet$\hss}}%
    \graphtemp=.5ex\advance\graphtemp by 0.111in
    \rlap{\kern 0.333in\lower\graphtemp\hbox to 0pt{\hss $\bullet$\hss}}%
    \special{pn 8}%
    \special{pa 778 111}%
    \special{pa 1000 496}%
    \special{pa 778 881}%
    \special{pa 333 881}%
    \special{pa 111 496}%
    \special{pa 333 111}%
    \special{pa 778 111}%
    \special{fp}%
    \special{pa 1000 496}%
    \special{pa 333 881}%
    \special{pa 333 111}%
    \special{pa 1000 496}%
    \special{fp}%
    \graphtemp=.5ex\advance\graphtemp by 0.111in
    \rlap{\kern 1.667in\lower\graphtemp\hbox to 0pt{\hss $\bullet$\hss}}%
    \graphtemp=.5ex\advance\graphtemp by 0.496in
    \rlap{\kern 1.889in\lower\graphtemp\hbox to 0pt{\hss $\bullet$\hss}}%
    \graphtemp=.5ex\advance\graphtemp by 0.881in
    \rlap{\kern 1.667in\lower\graphtemp\hbox to 0pt{\hss $\bullet$\hss}}%
    \graphtemp=.5ex\advance\graphtemp by 0.881in
    \rlap{\kern 1.222in\lower\graphtemp\hbox to 0pt{\hss $\bullet$\hss}}%
    \graphtemp=.5ex\advance\graphtemp by 0.496in
    \rlap{\kern 1.000in\lower\graphtemp\hbox to 0pt{\hss $\bullet$\hss}}%
    \graphtemp=.5ex\advance\graphtemp by 0.111in
    \rlap{\kern 1.222in\lower\graphtemp\hbox to 0pt{\hss $\bullet$\hss}}%
    \special{pa 1667 111}%
    \special{pa 1889 496}%
    \special{pa 1667 881}%
    \special{pa 1222 881}%
    \special{pa 1000 496}%
    \special{pa 1222 111}%
    \special{pa 1667 111}%
    \special{fp}%
    \special{pa 1667 111}%
    \special{pa 1667 881}%
    \special{pa 1000 496}%
    \special{pa 1667 111}%
    \special{fp}%
    \hbox{\vrule depth0.992in width0pt height 0pt}%
    \kern 2.000in
  }%
}%
}
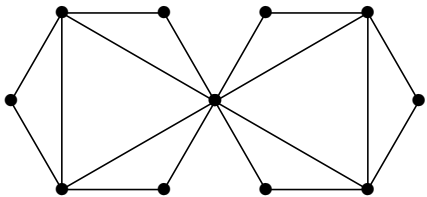
\ce{Fig.~1.  A graph with leafage 4 but no asteroidal 4-set}

\SH
{3. DOMINATORS AND A GENERAL UPPER BOUND}
As observed in the introduction, when studying $l(G)$ we may assume a bijection
between the maximal cliques of $G$ and the vertices of the host.  There is a
natural order on the maximal cliques.  If $\VEC vn1$ is a \peo\ of $G$, then the
reverse ordering $\VEC v1n$ constructs $G$ by iteratively adding a vertex
adjacent to a clique.  With respect to (the reverse of) a given \peo, we say
that a maximal clique is created whenever the vertex being added increases the
number of maximal cliques in the graph that has been built.  Each additional
vertex enlarges one maximal clique or creates one new maximal clique.  This
defines a {\it creation ordering} $\VEC Q1m$ of the maximal cliques of $G$ (with
respect to a given \peo).

\LM 4.
Let $\VEC Q1m$ be the creation ordering of the maximal cliques of a connected
chordal graph $G$ with respect to some \peo.  For each $i$ with $2\le i\le m$,
there is an index $j<i$ such that $Q_i\cap Q_k\esub Q_i\cap Q_j$ for all $k<i$.
We call $Q_j$ the {\it dominator} of $Q_i$ and write $j=\rho(i)$.
\PF
By induction on $n$, the order of $G$; the claim is vacuous for $n=1$.  For
$n>1$, apply the induction hypothesis to $G-v_n$, where $v_n$ is the first
vertex of the \peo.  When $v_n$ is added, it belongs to only one maximal clique.

Either $v_n$ enlarges a previous maximal clique, or it creates a new maximal
clique.  If $v_n$ enlarges a maximal clique $Q_i$, then the intersection
condition holds with the same choice of dominators as for $G-v_n$, because $v_n$
belongs to no other maximal clique.  If $v_n$ creates a new maximal clique
$Q_m$, then $N(v_n)$ is a proper subset of some earlier clique $Q_j$; set
$j=\rho(m)$ (this choice need not be unique).  Since $N(v_n)$ contains the
\intn\ of $Q_m$ with each other clique, the dominator condition holds.  \qed

\vskip -1pc
\LM 5.
Let $\VEC Q1m$ be the creation ordering of the maximal cliques of $G$ associated
with some \peo.  Let $T$ be the tree with vertices $\VEC q1m$ and edges $q_iq_j$
such that $j=\rho(i)$.  Assigning $q_i$ to all vertices of $Q_i$ yields a
subtree \rep\ of $G$ in $T$.  Furthermore, each clique of $G$ corresponding to a
leaf of $T$ contains a vertex appearing in no other maximal clique of $G$. 
\PF
By induction on $n$, the order of $G$; the claim is trivial for $n=1$.  For
$n>1$, apply the induction hypothesis to $G-v_n$, where $v_n$ is the first
vertex in the \peo.  If adding $v_n$ enlarges a maximal clique $Q_i$, then it
suffies to enlarge the \rep\ $f$ by setting $f(v_n)=q_i$.  If adding $v_n$
creates a new maximal clique $Q_m$ with dominator $Q_j$, then $T$ gains the
pendant edge $q_m q_j$ and is still a tree.  Leaves still contain unique
representatives, since we let $f(v_n)=q_m$.  Finally, since $Q_j$ contains
$N(v_n)$, the subtrees assigned to the vertices of $N(v_n)$ can be extended from
$q_j$ to include $q_m$.  Hence we obtain the desired subtree \rep\ of $G$.  \qed

To obtain an upper bound on $l(G)$, we consider a partial order associated with
$G$.  An {\it antichain} in a partial order is a set of pairwise incomparable
elements, and the {\it width} $w(P)$ of a partial order $P$ is the maximum size
of its antichains.  Dilworth's Theorem [4] states that the elements of a finite
partial order of width $k$ can be partitioned into $k$ chains (totally ordered
subsets).  Such a set of chains is a {\it Dilworth decomposition}.

For each simplicial neighborhood $R$ in a chordal graph $G$, we define the {\it
modified simplicial neighborhood} to be $R'=R-S(G)$.  Let $\bR'(G)$ denote the
set of modified simplicial neighborhoods in $G$.  For $v\in S(G)$, we also
define $N'(v)=N(v)-S(G)$.

\TH 4.
If $P(G)$ is the inclusion order on the collection $\bR'(G)$ of modified
simplicial neighborhoods in a chordal graph $G$, then $l(G)\le w(P(G))$.
\PF
Let $\VEC C1k$ be a Dilworth decomposition of $P(G)$, with $C_j$ consisting of
$\VECOP R{j1}{jr_j}\esub$.  For each chain $C_j$ we create a path
$P_j$ to be part of the host tree in a subtree \rep\ $f$.  If $u,v\in S(G)$ are
adjacent, then $N'(u)=N'(v)$.  For $R\in \bR'(G)$, let $m(R)$ be the number of
components of $S(G)$ for which $R$ is the common modified simplicial
neighborhood.  We assign to $R_{ji}$ exactly $m(R_{ji})$ consecutive vertices
on $P_j$, each of which will be the entire image for the vertices of one
component of $S(G)$.  We put the vertices for $R_{j1}$ at one end of the path
(the {\it small} end, pass through those for each $R_{ji}$ as $i$ increases, and
reach the vertices for $R_{jr_j}$ at the {\it big} end of the path.  Each vertex
of the subpath for $R_{ji}$ is assigned to each vertex in $R_{ji}$.  For each
$v'\in V(G')$, we have assigned to $v'$ a (possibly empty) subpath from the big
end of each $P_j$.

The problem is now to add host vertices to hook together the big ends of these
paths so that 1) the only leaves are the small ends of $\VEC P1k$, 2) the
maximal cliques not containing simplicial vertices of $G$ are represented, and
3) the vertices assigned to each vertex of $G'$ form a subtree of the resulting
tree.

The maximal cliques of $G$ that do not contain simplicial vertices of $G$ are
maximal cliques of $G'$.  Let $f'$ be the subtree \rep\ of $G'$ generated by
Lemma 5 in a host tree $T'$.  In $T' \cup \{\VEC P1k\}$ all maximal cliques of
$G$ are represented.  Each ``big end'' $R_{jr_j}$ is contained in (possibly with
equality) a maximal clique $Q$ of $G'$; add an edge from the big end of $P_j$ to
the vertex for $Q$ in $T'$.  (If $R_{jr_j}=Q$, then this edge can be
contracted.)  The result is a tree $T$.  Each $v'\in V(G')$ was assigned a
subtree in $T'$; if $v'$ received additional vertices in $P_j$, then these have
been attached to a vertex of $f'(v')$, so the image $f(v')$ is still a tree.

Finally, we show that each leaf of $T'$ receives an edge from some path $P_j$;
thus all leaves in the full tree $T$ are small ends of $\VEC P1k$.  Let $q$ be a
leaf of $T'$.  By Lemma 5, the clique of $G'$
corresponding to $q$ contains a vertex in no other maximal clique of $G'$.
Such a vertex $v'$ is simplicial in $G'$.  Since $G'$ contains no simplicial
vertices of $G$, $v'$ has a neighbor $v\in S(G)$.  Let $C_j$ be the chain in
$P(G)$ containing $N(v)$.  The big end of $P_j$ must be attached to $q$, because
$R_{jr_j}$ contains $v'$ and $q$ is the only vertex of $T'$ assigned to $v'$.
\qed

\SH
{4. CHARACTERIZATION OF MAXIMUM ASTEROIDAL SETS}
The sharpness of the upper bound in Theorem 4 depends on the structure of $G'$.
When $G'$ is a single clique, it has no cutset, and Lemma 6 below implies that
$a(G)=l(G)=w(P(G))$.  Among such graphs we find the $n$-vertex chordal graph
with maximum leafage.  Lemma 6 is a lower bound on $a(G)$ in terms of a subposet
of $P(G)$.  Given a chordal graph $G$, the {\it restricted simplicial
neighborhood poset} $P'(G)$ is obtained from $P(G)$ by deleting those modified
simplicial neighborhoods that contain cutsets of the derived graph $G'$.

\LM 6.
If $P'(G)$ is the restricted simplicial neighborhood poset of a chordal graph
$G$, then $G$ has an asteroidal set consisting of $w(P'(G))$ simplicial
vertices.
\PF
Let $\VEC R1r$ be a maximum antichain in $P'(G)$.  For each $R_i$, choose
$v_i\in S(G)$ with $N'(v_i)=R_i$.  Given distinct vertices $v_i,v_j,v_k$ in the
resulting set, we may choose $v_j'\in R_j-R_i$ and $v_k'\in R_k-R_i$ because
$\VEC R1r$ is an antichain.  By the definition of $P'(G)$, $G'-R_i$ is connected
and contains a $v_j',v_k'$-path $P$.  Appending $v_j$ and $v_k$ yields a
$v_j,v_k$-path avoiding $N(v_i)$.  \qed

\vskip -1pc
\TH 5.
The maximum leafage of a chordal graph with $n$ vertices is the maximum $k$
such that $k\le \CH{n-k}{\FL{(n-k)/2}}$.  In terms of $n$, the value is
$n- \lg n -\half \lg\lg n +O(1)$.
\PF
When $k$ satisfies the inequality, we construct a chordal $n$-vertex graph with
an asteroidal set of size $k$.  To a clique $Q$ of order $n-k$, we add $k$
simplicial vertices whose neighborhoods are distinct $\FL{(n-k)/2}$-subsets of
$G'=Q$.  By Lemma 6, $S(G)$ is an asteroidal set of size $k$, and hence
$l(G)\ge k$ by Lemma 1.

Let $G$ be an $n$-vertex chordal graph with maximum leafage.  Theorem 4 implies
that $l(G)\le w(P(G))$.  Let $m=\C{S(G)}$.  Since the elements of $P(G)$ are
distinct modified simplicial neighborhoods, $w(P(G))\le m$.  Since elements
$P(G)$ are subsets of an $n-m$-element set, $w(P(G))\le \CH{n-m}{\FL{(n-m)/2}}$.
Thus $w(P(G))\le \min\{m,\CH{n-m}{\FL{(n-m)/2}}\}$.  Since the second term
decreases with $m$, the width cannot exceed the maximum $k$ such that
$k\le\CH{n-k}{\FL{(n-k)/2}}$.  This proves the upper bound.  To obtain the
asymptotic maximum, we apply Stirling's approximation to
$k = \CH{n-k}{\FL{(n-k)/2}}$.  \qed

We can improve the lower bound resulting from Lemma 6 by applying Lemma 6 to
every connected induced subgraph of $G$.

\TH 6.
If $G$ is a chordal graph, then the maximum size of an \ast\ set in $G$ is
the maximum of $w(P'(H))$ over all connected induced subgraphs $H$ of $G$.
\PF
Every induced subgraph $H$ of a chordal graph $G$ is chordal.  If $H$ has an
$x,y$-path $P$ avoiding $N(z)$ for some $x,y,z\in V(H)$, then $P$ is also an
$x,y$-path avoiding $N(z)$ in $G$, since $G$ has no additional edges among
vertices of $H$.  Thus every asteroidal set in $H$ is an asteroidal set in
$G$, and the lower bound follows from Lemma 6.

For the upper bound we construct, from a maximum \ast\ set $A$ in $G$, a
connected induced subgraph $H$ of $G$ such that $\C A = w(P'(H))$.
For each triple $x,y,z \in A$, there is an $x,y$-path in $G$ avoiding $N(z)$.
We may assume that each such path is chordless, since shortening the path by
using chords still avoids $N(z)$.  Let $H$ be the subgraph of $G$ induced by
the vertices in the union of all these chordless paths.  It suffices to show
that 1) $S(H)=A$, 2) $P(H)$ is an antichain of size $\C A$, and 3) $P'(H)=P(H)$.

Since $H$ is an induced subgraph of $G$, $N_H(x)\esub N_G(x)$, so the vertices
of $A$ are simplicial in $H$.  All members of $V(H)-A$ are internal vertices
of chordless paths between vertices of $A$.  Such a vertex cannot be simplicial
in $H$, because its neighbors on such a path are not adjacent; this proves (1).
If $N_H(x)\esub N_H(y)$ for some $x,y\in A$, then every $x,z$-path intersects
$N(y)$; since $A$ is asteroidal, this proves (2).

It remains only to show that $P'(H)=P(H)$.  First note that $N_H'(x)=N_H(x)$ for
$x\in S(H)$ since $S(H)=A$ is an independent set.  Thus it suffices to show that
$N_H(x)$ does not contain a cutset of $H'$ for $x\in A$.  Suppose that $Q$ is a
minimal cutset of $H'$ contained in $N_H(x)$.  Since $Q$ is a minimal cutset of
a chordal graph $H'$, $Q$ induces a clique, and every component of $H'-Q$
contains a simplicial vertex of $H'$.  Each such vertex is a neighbor of a
simplicial vertex of $H$.  This yields vertices $y,z\in A$ such that every
$y,z$-path in $H$ intersects $Q$, which contradicts $A$ being an asteroidal set.
\qed

For the graph $G$ in Fig.~2, the four simplicial vertices of degree 2 establish
$a(G)\ge4$.  Equality holds; these vertices and their neighborhoods together
induce the only connected induced subgraph $H$ such that $w(P'(H))=4$.  Although
$w(P'(H))$ is computable from $H$ in polynomial time, the difficulty of finding
$H$ makes the complexity of computing $a(G)$ unclear.  For example, if $F$ is
formed by adding a pendant edge to each leaf of the graph in Fig.~2, then
$w(P'(F))=2$ and $a(F)=4$, but $F$ has no simplicial vertex $x$ such that
$w(P'(F-x))>2$.  Thus $H$ cannot be found greedily.

\vskip -1pc
\gpic{
\expandafter\ifx\csname graph\endcsname\relax \csname newbox\endcsname\graph\fi
\expandafter\ifx\csname graphtemp\endcsname\relax \csname newdimen\endcsname\graphtemp\fi
\setbox\graph=\vtop{\vskip 0pt\hbox{%
    \graphtemp=.5ex\advance\graphtemp by 0.865in
    \rlap{\kern 0.865in\lower\graphtemp\hbox to 0pt{\hss $\bullet$\hss}}%
    \graphtemp=.5ex\advance\graphtemp by 0.865in
    \rlap{\kern 1.635in\lower\graphtemp\hbox to 0pt{\hss $\bullet$\hss}}%
    \graphtemp=.5ex\advance\graphtemp by 0.096in
    \rlap{\kern 1.635in\lower\graphtemp\hbox to 0pt{\hss $\bullet$\hss}}%
    \graphtemp=.5ex\advance\graphtemp by 0.096in
    \rlap{\kern 0.865in\lower\graphtemp\hbox to 0pt{\hss $\bullet$\hss}}%
    \graphtemp=.5ex\advance\graphtemp by 0.865in
    \rlap{\kern 1.250in\lower\graphtemp\hbox to 0pt{\hss $\bullet$\hss}}%
    \graphtemp=.5ex\advance\graphtemp by 0.481in
    \rlap{\kern 1.635in\lower\graphtemp\hbox to 0pt{\hss $\bullet$\hss}}%
    \graphtemp=.5ex\advance\graphtemp by 0.096in
    \rlap{\kern 1.250in\lower\graphtemp\hbox to 0pt{\hss $\bullet$\hss}}%
    \graphtemp=.5ex\advance\graphtemp by 0.481in
    \rlap{\kern 0.865in\lower\graphtemp\hbox to 0pt{\hss $\bullet$\hss}}%
    \special{pn 8}%
    \special{pa 865 865}%
    \special{pa 1635 865}%
    \special{pa 1635 96}%
    \special{pa 865 96}%
    \special{pa 865 865}%
    \special{fp}%
    \special{pa 1250 865}%
    \special{pa 1635 481}%
    \special{pa 1250 96}%
    \special{pa 865 481}%
    \special{pa 1250 865}%
    \special{fp}%
    \graphtemp=.5ex\advance\graphtemp by 0.481in
    \rlap{\kern 0.096in\lower\graphtemp\hbox to 0pt{\hss $\bullet$\hss}}%
    \graphtemp=.5ex\advance\graphtemp by 0.481in
    \rlap{\kern 0.481in\lower\graphtemp\hbox to 0pt{\hss $\bullet$\hss}}%
    \graphtemp=.5ex\advance\graphtemp by 0.481in
    \rlap{\kern 2.019in\lower\graphtemp\hbox to 0pt{\hss $\bullet$\hss}}%
    \graphtemp=.5ex\advance\graphtemp by 0.481in
    \rlap{\kern 2.404in\lower\graphtemp\hbox to 0pt{\hss $\bullet$\hss}}%
    \special{pa 96 481}%
    \special{pa 865 481}%
    \special{fp}%
    \special{pa 1250 865}%
    \special{pa 1250 96}%
    \special{fp}%
    \special{pa 1635 481}%
    \special{pa 2404 481}%
    \special{fp}%
    \hbox{\vrule depth0.962in width0pt height 0pt}%
    \kern 2.500in
  }%
}%
}
\ce{Fig.~2.  A graph with $a(G)>w(P'(G))$.}

\skipit{
It is tempting to try to use this theorem to define an algorithm to seek a
maximum \ast\ set in $G$.  Given $x\in S(G)$, we can test whether
$w(P'(G-x))\ge w(P'(G))$.  If it is, we delete $x$ and iterate.  We stop
and select the asteroidal set generated by an antichain in $P'$ of the remaining
graph $H$ when $w(P'(H-x))\le w(P'(H))$ for all $x\in S(H)$.  Unfortunately,
this greedy procedure does not work.  
}

In the computation of $a(G)$, it may be useful to confine the maximum asteroidal
set to the set $S(G)$ of simplicial vertices.

\LM 7.
If $G$ is a chordal graph, then $S(G)$ contains a maximum \ast\ set of $G$.
\PF
Let $A$ be a maximum asteroidal set of $G$.  If $x\in A$ is not simplicial
in $G$, then $x$ has nonadjacent neighbors.  Hence $x$ belongs to a minimal
cutset of $G$.  Since $G$ is chordal, this cutset induces a clique $Q$.
As in the proof of Theorem 6, $A-x$ is confined to one component of $G-Q$.
Any other component of $G-Q$ contains a simplicial vertex of $G$ that is
not in $Q$.  Choose such a vertex $w$ to replace $x$ in $A$.  To verify that
$A-x+w$ is an \ast\ set, it suffices to note that if $y,z\in A$ and $P$ is
a $y,x$-path avoiding $N(z)$, then $P$ can be extended from $x$ to $w$ to
obtain a $y,w$-path avoiding $N(z)$.  \qed

\SH
{5. LEAFAGE WHEN $G'$ HAS (AT MOST) TWO CLIQUES}
We have observed that $l(G)=a(G)=w(P(G))$ when $G'$ is a clique.  Computing
leafage in general means determining the savings $w(P(G))-l(G)$ in constructing
an optimal subtree \rep.  By studying antichains in more detail, we will be
able to compute leafage when $G'$ is the union of two maximal cliques.

An {\it ideal} of a poset $P$ is a subposet $Q$ such that all elements below
elements of $Q$ also belong to $Q$.  The inclusion order on the set of ideals of
$P$ is a lattice, because the union and intersection of two ideals is an ideal.
Every antichain in $P$ is associated naturally with the ideal for which it forms
the set of maximal elements.  The ordering of the antichains that corresponds to
the inclusion ordering on the ideals has $A\le B$ when for every $x\in A$ there
exists $y\in B$ such that $x\le y$.

Dilworth [5] proved that the set of maximum-sized antichains also forms a
lattice under this ordering.
\skipit{
If $A$ and $B$ are maximum antichains, then the set of
minimal elements in $A\cup B$ and the set of maximal elements in $A\cup B$ are
the antichains corresponding to the intersection and union of the ideals
corresponding to $A$ and $B$.  Since they both contain $A\cap B$, they both also
have size $w(P)$.
}
When $P$ is finite, the lattice of maximum-sized antichains is finite and has a
unique minimal element, which we denote by $A_P$.  In an arbitrary poset, a
maximum antichain can be found quickly using network flow methods or bipartite
matching.  It is well-known that in the bipartite graph $\GA(P)$ with partite
sets $\{x^-\st x\in P\}$ and $\{x^+\st x\in P\}$ and edge set
$\{x^-y^+\st x<y {\rm\ in\ } P\}$, the maximum size of a matching is
$\C P -w(P)$.  Using this, we can find $A_P$.

\LM 8.
In a finite poset $P$, the unique minimal maximum antichain $A_P$ can be
found in polynomial time.
\PF
\skipit{
Let $\bA$ be the lattice of maximum antichains in $P$.  Using bipartite
matching on $\GA(P)$, we find some $A\in\bA$.  Let $A^-$ be the ideal of $P$ for
which $A$ is the antichain of maximal elements.  No element of $P$ above any
element of $A$ can belong to $A_P$, since each element of $A_P$ is below some
element of $A$; hence $A_P\esub A^-$.
|
If $A^-$ contains any antichain $A'$ of size $w(P)$ other than $A$, then it is
below $A$ in $\bA$, since each element of $A^-$ is below some element of $A$.
If there is any such antichain, it must omit some element of $A$.  Hence we run
bipartite matching on $\GA(A^- - x)$ for each $x \in A$.  If some such run finds
an antichain of size $w(P)$, let it be the new $A$ and repeat this procedure
with the new $A$ and new $A^-$.  If all such runs find no antichain of size
$w(P)$, then $A=A_P$.  Since $\C{A^-}$ decreases with each iteration, this
algorithm terminates by finding $A_P$ in at most $w(P)\C P$ uses of the
bipartite matching algorithm.
{\it Alternative proof.}
It is also possible to find $A_P$ with only one use of bipartite matching but
more post-processing, because the edges in a maximum matching of $\GA(P)$ 
}
Using a bipartite matching algorithm, we find a maximum matching of $\GA(P)$.
These edges link the elements of $P$ into chains that form a Dilworth
decomposition $\VEC C1{w(P)}$ of $P$.  Every maximum antichain uses one element
from each $C_i$.  If the minimal elements of the chains form an antichain, then
this is $A_P$.  Otherwise, there is some relation $x > y$ among these elements.
Since $y$ is below every element of the chain containing $x$, no maximum
antichain contains $y$.  Since $y$ belongs to no maximum antichain in $P$, we
have $w(P-y)=w(P)$, and thus the same chain partition with $y$ deleted forms a
Dilworth decomposition of $P-y$.   Iterating this procedure with the minimal
remaining elements eventually produces an antichain of size $w(P)$.
Since no maximum antichains of $P$ were destroyed, this set is $A_P$.  \qed

\LM 9.
Let $\bC$ be a Dilworth decomposition of $P$.  If $x\in A_P$ and $C'$ is the
chain consisting of $x$ and the elements above $x$ on its chain in $\bC$, then
$w(P-C')=w(P)-1$.  Equivalently, $P$ has a Dilworth decomposition using
$C'$ as one chain.
\PF
Let $k=w(P)$.  Since $A_P-x$ is an antichain in $P-C'$, it suffices to show that
$P-C'$ has no antichain of size $k$.  Let $C$ be the chain containing $x$ in
$\bC$.  If $P-C'$ has an antichain $A$ of size $k$, then $A$ must contain an
element $y$ below $x$ on $C$, since $w(P-C)<k$.  Also $A$ is an antichain of
size $k$ in $P$, and hence $A_P<A$.  In particular, $x$ is less than some member
of $A$, which creates a forbidden relation between $y$ and another member of
$A$.  \qed


\vskip -1pc
\LM {10}.
Let $Q_1,Q_2$ be distinct maximal cliques of the derived graph $G'$ of a chordal
graph $G$.  Let $P=P(G)$, and define induced subposets
$P_1 = \{x\in P\st Q_1\cap Q_2\subset x\esub Q_1\}$,
$P_2 = \{x\in P\st Q_1\cap Q_2\subset x\esub Q_2\}$, and $\Pb = P_1 \cup P_2$.
Suppose that $A_P$ has elements $a_1\in P_1$ and $a_2\in P_2$.
If $w(P)\ge w(P-\Pb)+2$, then $P(G)$ has a Dilworth decomposition in which
$a_1,a_2$ occur as minimal elements of chains.
\PF
With these hypotheses, Lemma 9 implies that $P$ has a Dilworth decomposition
containing a chain $C_1$ with $a_1$ as its minimal element.  Let $P^*=P-C_1$,
and let $A^*=A_P-\{a_1\}$.  Since $A^*$ is a maximum antichain of $P^*$, we have
$A_{P^*}\le A^*$ in the lattice of maximum antichains of $P^*$.

Let $\GA$ be the bipartite graph with partite sets $A_{P^*}$ and $A^*$ (each
of size $w(P)-1$) defined by putting $x\in A_{P^*}$ adjacent to $y\in A^*$ if
$x\le y$ in $P$.  If $\GA$ has a vertex cover with fewer than $w(P^*)$ vertices,
then the vertices not in the cover form an independent set, and this independent
set forms an antichain of size exceeding $w(P^*)$ in $P^*$.  Hence there is no
such cover, and by the K\"onig-Egerv\'ary Theorem $\GA$ has a complete matching
$M$.

Let $B=A_P\cap P_2$; note that $B \esub A^*$.  For $b\in B$, let $b'$ be the
element matched to $b$ in $M$, and let $B'=\{b'\st b\in B\}$.  We claim that
$B'\esub P_2$.  For $b\in B$, we have $b'\le b$ in $P$, so
$b'\esub b\in B\esub P_2$.  Thus $b'\esub Q_2$.  If $b'\esub Q_1\cap Q_2$, then
$b'$ is properly contained in every element of $\Pb$.  Since $b'\in A_{P^*}$,
this implies that $A_{P^*}$ is disjoint from $\Pb$, and we have
$w(P-\Pb)\ge \C{A^*}=w(P)-1$.  We assumed that $w(P)\ge w(P-\Pb)+2$, so $b'$
must intersect $Q_2-Q_1$.  This implies that $b'\in P_2$.

With $B'\esub P_2$, we next claim that $B'=B$.  A set not in $P_2$ that contains
an element of $P_2$ contains $Q_2$.  Since $a_2\in A_P\cap P_2$, no such element
is in the antichain $A_P$.  Thus $b'$ is not below any element of $A_P-B$.  Also
$b'$ is not above any element of $A_P-B$, because $b'\esub b$ and $A_P$ is an
antichain.  Hence $A'=(A_P-B)\cup B'$ is an antichain in $P$, with size $w(P)$
since $\C B=\C{B'}$.  Since $B'\le B$, we have $A'\le A_P$ in the lattice.
Since $A_P$ is the minimal element of the lattice, this yields $A'=A_P$ and
$B'=B$.

Since $B'=B$, we have proved that all the elements of $A_P$ belonging to $P_2$
also belong to $A_{P^*}$.  In particular, $a_2\in A_{P^*}$.  Applying Lemma 9
again, we obtain a Dilworth decomposition of $P^*$ containing a chain $C_2$ with
$a_2$ as its minimal element.  Together with $C_1$, this is the desired Dilworth
decomposition of $P$.  \qed

Recall that $P'(G)$ is the subposet of $P(G)$ obtained by discarding the
modified simplicial neighborhoods that contain cutsets of $G'$.  When the
modified simplicial neighborhoods contained in a particular maximal clique $Q_i$
of $G'$ form a chain under inclusion, we say that $Q_i$ is {\it degenerate}. 

\TH 7.
Let $G$ be a chordal graph such that $G'$ is the union of two maximal
cliques $Q_1,Q_2$.  With $P=P(G)$, $P'=P'(G)$, and $Q = Q_1\cap Q_2$, $l(G)$
has the following value:
\item{1)} $l(G)=w(P)-\al$ if $w(P')\le w(P)-2$ and $\al$ of the cliques
$Q_1,Q_2$ are nondegenerate and contain distinct elements of $A_{P}-P'$.
\item{2)} $l(G) = w(P)$ if $w(P')=w(P)-1$ and every element of
$A_{P}-P'$ belongs to a degenerat $Q_i$.
\item{3)} $l(G)=w(P')$ otherwise.
\PF
Lemma 6 and Theorem 4 imply that $w(P')\le l(G)\le w(P)$, so we may assume that
$w(P')<w(P)$.  The only minimal cutset of $G'$ is $Q$, so the elements of $P$
discarded to form $P'$ are those containing $Q$.  Also, every element of
$P$ is contained in $Q_1$ or in $Q_2$.

We first consider improvements to the upper bound.  If $w(P')\le w(P)-1$, then
$A_P$ has an element $x$ outside $P'$, which means $Q\esub x$.  We may assume
that $x \esub Q_1$; we can improve the upper bound if $Q_1$ is nondegenerate.
By Lemma 9, we can find a Dilworth decomposition of $P$ such that one chain has
$x$ as its bottom element.  Each chain consists of subsets of $Q_1$ or consists
of subsets of $Q_2$.  Using the construction in the proof of Theorem 4, we
produce two host subtrees, with a total of $w(P)$ leaves together (if $Q_2$ is
also nondegenerate), one of which has vertices for simplicial neighborhoods
contained in $Q_1$, the other for $Q_2$.  Furthermore, $x$ appears at a leaf in
the tree for $Q_1$.  Since $x$ contains $Q$, we can add an edge between the
vertex for $x$ and the vertex representing $Q_2$ in the other tree.  This yields
a subtree \rep\ of $G$ with $w(P)-1$ leaves.  (If $Q_2$ is degenerate, then the
initial pair of subtrees has a total of $w(P)+1$ leaves, but the tree for
subsets of $Q_2$ is a path with the vertex for $Q_2$ as a leaf, and the added
edge eliminates two leaves.)

If $w(P')\le w(P)-2$ and $Q_1,Q_2$ contain distinct elements of $A_P-P'$, then
we may be able to save another leaf.  The hypotheses of Lemma 10 hold (with
$P'=P-\Pb$), and Lemma 10 guarantees a Dilworth decomposition of $P$ having
chains with bottom elements $x,y$ such that $x,y\in P'$, $Q\esub x\esub Q_1$,
and $Q\esub y\esub Q_2$.  As above, we use the construction in the proof of
Theorem 4 to produce representations in two host subtrees, one having $x$ at a
leaf and the other having $y$ at a leaf.  The two subtrees have a total of
$w(P)+2-\al$ leaves (each of $Q_1,Q_2$ that is degenerate increases the number
of leaves in the initial pair of trees by one).  By adding the edge $xy$, we
produce a subtree \rep\ of $G$ with $w(P)-\al$ leaves.

We next prove that $w(P)-\al$ is a lower bound when $w(P')\le w(P)-2$.  The
latter inequality implies that no element of $A_P$ is contained in $Q$.  Hence
$m_1+m_2=w(P)$, where $A_P$ consists of $m_1$ subsets of $Q_1$ and $m_2$ subsets
of $Q_2$.  Suppose first that $m_1,m_2 >0$.  By the argument in Lemma 6, the
simplicial vertices of $G$ corresponding to the elements of $A_P$ contained in
$Q_i$ yield an asteroidal set $U_i$ of size $m_i$.  Consider an optimal
representation of $G$.  Since the
host has a vertex for each maximal clique of $G$ and simplicial vertices of $G$
appear in exactly one maximal clique, each vertex of $U_i$ is assigned one host
vertex.  Let $T_i$ be the subtree of the host consisting of all paths between
vertices represesenting $U_i$.  Since $U_i$ is an asteroidal set, Theorem 1
implies that the corresponding subtrees form an asteroidal collection.  Since
these trees are single vertices, the leaves of $T_i$ are precisely the $m_i$
vertices assigned to elements of $U_i$.

If the vertex assigned to $u\in U_i$ lies between the vertices assigned to
$v,w\in U_{3-i}$, then every $v,w$-path in $G$ contains a neighbor of $u$.
On the other hand, the vertices $v,w$ have neighbors $v',w' \in Q_{3-i}-Q_i$,
and the path $v,v',w',w$ avoids $N(u)$.  This contradiction implies that
the trees $T_1,T_2$ are disjoint.  We now have a total of $w(P)$ leaves in two
disjoint subtrees, except that this total increases by one for each $Q_i$
that is degenerate or contains no element of $A_P-P'$.  The host contains a
unique path between these subtrees, which reduces the number of leaves by at
most two.  Hence the host has a subtree with $w(P)-\al$ leaves, and
$l(G)\ge w(P)-\al$.

When $w(P')=w(P)-1$ and every element of $A_P-P'$ belongs to a degenerate $Q_i$,
a similar argument shows that $l(G)=w(P)-1$.  \qed

\skipit{
When $m_2=0$ (similarly for $m_1=0$), we have $Q\subset x\esub Q_1$ for all
$x\in A_P-P'$.  The argument above yields $T_1$ as before with a
total of $w(P)$ leaves, but there is no \rep\ of $Q_2$ using only $m_2$ leaves
that has a set containing $Q$ at a leaf (because the subtrees corresponding
to $U_2$ form an asteroidal collection), and hence the edge added to join
the two subtrees into a single tree by uniting the subtrees for $Q$ in each
$T_i$ can save only one leaf.
\qed
}

Perhaps these ideas can be combined with the ``dominator tree'' to obtain a
polynomial-time algorithm in general.  The condition that $G'$ has at most two
cliques is a recognizable, since the algorithm of Rose, Tarjan, and Leuker [21]
finds a \peo\ (if one exists) in time linear in the number of vertices plus
edges, and from this the maximal cliques and simplicial vertices are available.
An algorithm for recognizing chordal graphs with leafage at most 3 can be
obtained from the material in [18].

\skipit{Meanwhile, we can use existing algorithms to recognize the graphs with
leafage at most 3.
\TH 8.
The leafage of $G$ is at most three if and only if $G$ is the union of graphs
$G_1 , G_2$, with $G_1 \cap G_2$ being a clique with vertex set $S$, such that
each of $G_1,G_2$ is an interval graph and at least one of $G_1,G_2$ has an
interval \rep\ having the interval for all vertices of $S$ start before any
other interval starts.
\PF
If $l(G)=2$, the condition holds vacuously.  Suppose $l(G)=3$.  Then the host
tree $T$ in an optimal subtree \rep\ of $G$ has exactly one branch point $q$.
The vertices whose subtrees contain both endpoints of an edge incident to $q$
induce a clique $S$ that separates $G$.  The subhosts of $T$ obtained by 
deleting this edge provide the desired \rep s for the subgraphs $G_1,G_2$.
    Conversely, if $G$ has the decomposition described, then one end of the 
interval \rep\ for one of the subgraphs can be attached to the interval
\rep\ for the other subgraph to obtain a subtree \rep\ with at most three
leaves.  \qed
    This leads to a polynomial-time recognition algorithm because an $n$-vertex
chordal graph has fewer than $n$ minimal separators, which can be found along
with the maximal cliques while testing for the existence of a \peo\ in time
$O(e(G))$.  For each of the subgraphs thus formed, we check whether it is
an interval graph and whether it has an interval \rep\ in which the vertices of
$S$ appears at the beginning.  To test this property in $H$, we can append
a vertex $y$ with neighborhood $S\esub V(H)$ and another vertex $x$ with
neighborhood $y$ and test whether the new graph $H'$ is an interval graph.
Note also that if $G-S$ has more than two components, it does not matter
how the components are grouped to form $G_1$ and $G_2$.
}

\SH
{6. PROPER LEAFAGE AND EXTREME POINTS}

We now consider proper leafage.  The graphs with proper leafage 2 are the proper
interval graphs, which Roberts [19] proved are precisely the unit interval
graphs.  Less well-known is a structural \charz\ proved by Roberts;
we generalize this to obtain bounds on proper leafage.

The {\it closed neighborhood} of a vertex $a$ is the set $N[a]=N(a)\cup\{a\}$.
Vertices with the same closed neighborhood are {\it equivalent} in $G$; this
defines an equivalence relation on $V(G)$.  Each equivalence class induces a
clique, and vertices $x,y$ from distinct classes are adjacent if and only if
every vertex equivalent to $x$ is adjacent to every vertex equivalent to $y$.
The {\it reduction} $G^*$ of a graph $G$ is the subgraph induced by selecting
one vertex from each equivalence class.  A graph is {\it reduced} if its has no
pair of distinct equivalent vertices.

A vertex $a\in V(G)$ is an {\it extreme point} (EP) in $G$ if (1) $a$ is
simplicial and (2) every pair of vertices in $N(a)$ that are not equivalent to
$a$ have a common neighbor outside $N[a]$.  The simplicial vertices in $K_4-e$
are extreme points, but $K_4-e$ is not reduced since the three-valent vertices
have the same closed neighborhood.  Roberts used extreme points and reduced
subgraphs to characterize proper interval graphs.

\TH 8.
(Roberts [19])
A graph $G$ is a proper interval graph if and only if every connected reduced
induced subgraph of $G$ has at most two extreme points.  \qed

Roberts also proved that this statement holds when ``extreme points'' is
replaced by ``modified extreme points''.  If $H$ is a component of $G-N[a]$,
let $\partial H = \{x\in N[a]\st x$ has a neighbor in $H\}$.  We say that
a simplicial vertex $a$ is a {\it modified extreme point} (MEP) if
$\partial H_1 = \partial H_2$ for every pair of components $H_1,H_2$ of
$G-N[a]$.

\RK 1.
If $a$ is an MEP in a connected chordal graph $G$ that is not a clique, and $S$
is the set vertices of $N[a]$ not equivalent to $a$, then $S$ is a minimal
cutset of $G$.
\PF
Every vertex of $S$ has a neighbor in some component of $G-N[a]$.  If $a$ is an
MEP, then the vertices of $S$ with neighbors in each component of $G-N[a]$ are
the same.  Hence deleting any proper subset of $S$ does not separate $G$.  \qed

As observed by Gavril [8], every chordal graph has a subtree \rep\ by a
proper family of subtrees in a host tree, meaning that no subtree properly 
contains another.  Roberts' \charz\ works because modified extreme points force
leaves in a proper subtree \rep much as asteroidal sets force leaves in a
subtree \rep.  This does not hold for extreme points.  Roberts observed that the
concepts of EP and MEP are equivalent on claw-free connected chordal graphs.
The following structural lemma implies that MEP's in chordal graphs are EP's,
but we will see that the converse need not hold.  

\LM {12}.
If $G$ is a chordal graph, and $S$ is a minimal cutset of $G$, then
every component of $G-S$ contains a vertex adjacent to every vertex of $S$.
\PF
Let $H$ be a component of $G-S$, and let $x$ be a vertex of $H$ having the
maximum number of neighbors in $S$.  If there exists $v\in S-N(x)$, then
$v$ must have a neighbor in $V(H)$, since $S$ is a minimal cutset of $G$.
Choose $y\in N(v)\cap V(H)$ with minimal distance from $x$ in $H$, and let $P$
be a shortest $x,y$-path in $H$.  Since $G$ is chordal, $S$ induces a clique,
and hence every vertex of $S\cap N(x)$ completes a cycle with $P$ and $v$.
This leads to a chordless cycle unless $S\cap N(x) \subset N(y)$, which
contradicts the choice of $x$.  \qed

\vskip -1pc
\CO 2.
If $G$ is a connected chordal graph other than a clique, then every MEP in $G$
is an EP in $G$.
\PF
Suppose $a$ is an MEP of a $G$; note that both MEP's and EP's must be
simplicial.  By Remark 1, the set $S$ of vertices of $N[a]$ not equivalent
to $a$ is a minimal cutset.  By Lemma 12, every component of $G-S$ has a
vertex adjacent to all of $S$.  Hence any two vertices of $S$ have a common
neighbor outside $S$, and $a$ is an EP.  \qed



We next explore the role of MEP's in proper subtree \rep s.

\LM {13}.
Let $G$ be a non-clique chordal graph with a proper subtree \rep\ $f$ in a
host tree $T$.  If $v$ is an MEP of $G$ and $f(v)$ does not contain a leaf of
$T$, then $T-f(v)$ has a component disjoint from all subtrees representing
vertices outside $N[v]$.
\PF
Let $H$ be a component of $G-N[v]$.  The subtrees for $V(H)$ are confined to
a single component of $T-f(v)$, since their union is connected and shares no
vertex with $f(v)$.  If $w$ is a neighbor of $v$ not equivalent to $v$, then
$w$ has a neighbor in $H$, by Lemma 12.  Hence $f(w)$ contains the edge between
$f(v)$ and the component of $T-f(v)$ containing the subtrees for $V(H)$.
If each component of $T-f(v)$ contains a subtree for some vertex outside
$N[v]$, then $f(w)$ properly contains $f(v)$.  \qed

\vskip -1pc
\TH 9.
If $G$ is a connected non-clique chordal, and $H$ is a connected reduced induced
subgraph of $G$, then $l^*(G)$ is at least the number of MEP's in $H$.
\PF
A proper subtree \rep\ of $G$ must contain a proper subtree \rep\ of $H$.  It
suffices to show that every proper subtree \rep\ $f$ of $H$ has a distinct leaf
for each MEP of $H$.  Let $v$ be an MEP of $H$.  If $f(v)$ contains a leaf of
$T$, associate this leaf with $v$.  For example, if $T-f(v)$ is connected, then
there is only one edge from $f(v)$ to the rest of $T$, and $f(v)$ contains a
leaf of $T$.  If $T-f(v)$ has more than one component and $f(v)$ contains no
leaf, then Lemma 13 yields a component $T(v)$ of $T-f(v)$ that is disjoint from
the trees associated with the vertices of $G-N[v]$.  In this case, assign a leaf
of $T(v)$ to $v$.

It suffices to show that the leaves associated with distinct MEP's are distinct.
In a reduced graph, the simplicial vertices form an independent set.  Hence
MEP's $u,v$ are non-adjacent, which implies that $f(u),f(v)$ are disjoint.
Furthermore, $f(u)$ also cannot intersect a component $T(v)$ of $T-f(v)$ whose
vertices belong to subtrees in $f$ only for neighbors of $v$.  These two
statements imply that no leaf of $T$ belonging to $f(u)$ is associated with $v$.
Finally, suppose that neither $f(u)$ nor $f(v)$ contains a leaf of $T$.  In
addition to $f(u)\cap f(v)=\nul$, we have observed that $f(u)\cap T(v) =\nul$
and $f(v)\cap T(u) =\nul$.  Since the subtrees $T(u)$ and $T(v)$ are obtained
from $T$ by deleting an edge incident to $f(u)$ and to $f(v)$, respectively, we
conclude that they cannot have a common leaf.  \qed

\vskip -1pc
\LM {14}.
Every connected reduced chordal graph $G$ has an optimal proper subtree \rep\ in
which each MEP is assigned a leaf of the host tree that is assigned to no other
vertex.
\PF
Let $v$ be an MEP in $G$.  If $f(v)$ contains a leaf, we may extend $f(v)$ by
adding a new neighbor of the leaf and the leaf has the desired property.  If
$f(v)$ does not contain a leaf, then there is a component $T(v)$ of $T-f(v)$
whose vertices are assigned only to neighbors of $v$.  Since $G$ is reduced,
each neighbor of $v$ has neighbors outside $N[v]$; hence its subtree extends to
another component of $T-f(v)$.  We may therefore extend $f(v)$ to a leaf of
$T(v)$ and to one added vertex beyond it while maintaining a proper subtree
\rep.  \qed

To obtain the best lower bound from Lemma 14, we may need to consider
proper induced subgraphs.  The graph $P_5\join K_1$ has two simplicial
vertices, two MEP's, and two EP's, but it has an induced $\Koe$ (which
has three MEP's) and hence is not a proper interval graph.
Nevertheless, MEP's provide the right answer for block graphs (compare this
result with Theorem 3 and Corollary 1).
Note the contrast between this result and the fact that the leafage of a tree
$G$ is the number of leaves in $G'$.

\TH {10}.
The proper leafage of a block graph that is not a clique is the number of leaf
blocks.  In particular, the proper leafage of a tree is the number of leaves.
(The proper leafage of a clique is 2.)
\PF
In the reduction $G^*$ of a block graph $G$, each leaf block becomes an edge 
containing a simplicial vertex $v$.  This vertex $v$ is an MEP in $G^*$, because
$N(v)$ is a single vertex, and hence each component obtained by deleting $N[v]$
has that vertex as its neighborhood in $N[v]$.  This proves the lower bound; for
the upper bound, we construct a \rep.

If $G$ is a clique, we form a proper subtree \rep\ for $G$ by using a collection
of pairwise \intg\ subpaths of a path, with the initial vertices appearing
before all the terminal vertices and in the same order as the terminal vertices.
We may choose the two vertices represented at the leaves arbitrarily.

Let $m(G)$ denote the number of leaf blocks in $G$.  By induction on the number
of blocks in $G$, we build a proper subtree \rep\ with $\max\{m(G),2\}$ leaves,
having an arbitrary simplicial vertex from each leaf block appearing at the leaf
corresponding to that block (two such vertices if $G$ is a clique).  We have
verified this when $G$ is a clique.

If $G$ is choose a leaf block $B$, and let $v$ be the cut-vertex of $G$ in $B$.
Let $f'$ be the proper representation of $B$ in a path that has $v$ at one
leaf and the arbitrarily specified simplical vertex in $B$ and the other leaf.
If $v$ belongs to at least two blocks other than $B$, then $m(G-(B-v))=m(G)-1$.
The induction hypothesis provides a representation $f$ of $G-(B-v)$; we complete
this by adding an edge from the leaf assigned to $v$ in $f'$ to a vertex
assigned to $v$ in $f$.  If $v$ belongs to only one block other than $B$, then
this is a leaf block in $G-(B-v)$, and $m(G-(B-v))=m(G)$.  In this case the
representation $f$ of $G-(B-v))$ assigns a leaf to $v$, and we can add an edge
from it to the leaf assigned to $v$ in $f'$ without increasing the number of
leaves.  \qed

If $v$ is a simplicial vertex of a $k$-tree, then $N(v)$ is a minimal cutset,
and Lemma 12 then implies that $v$ is an MEP.  Hence the proper leafage of
a $k$-tree is at least the number of simplicial vertices.  The 2-tree
$P_5\join K_1$ shows that this bound need not be sharp.


%

\LM {15}.
Let $v$ be a cut vertex of a connected chordal graph $G$ such that $G-v$ has two
components.  If $G_1,G_2$ are the subgraphs obtained by deleting the vertices of
one of these components, then $l^*(G)\ge l^*(G_1)+l^*(G_2)-2$.  If $G_1$ or
$G_2$ has proper leafage 2 and is not a clique, then
$l^*(G)\ge l^*(G_1)+l^*(G_2)-1$.

\PF
Let $f$ be an optimal proper subtree \rep\ of $G$.  Since each $G_i$ is a
connected subgraph, the union of the subtrees assigned to vertices of $G_i$ is
a subtree of the host; call these $T_1,T_2$.  Furthermore, vertices of
$T_1\cap T_2$ can only be assigned to $v$.  If $v$ appears alone at a leaf, then
we could delete that leaf without losing the property of proper \rep.  Thus
we may assume that each leaf of $T$ is a leaf of exactly one of $T_1$ or $T_2$.
If we delete vertices assigned to any vertex outside $G_i$, we
still have a subtree \rep\ of $G_i$; we can guarantee that it is a proper \rep\
by growing a leaf from a vertex assigned to $v$.  Hence we have proved that
$l^*(G_1)+l^*(G_2)\le l^*(G)+2$.

If $l^*(G_2)=2$ and $G_2$ is not a clique, then we can obtain a
proper subtree \rep\ of $G_2$ with two leaves instead of growing an extra 
leaf from $T_2$ in the construction above.  Hence in this case we can improve
the bound to $l^*(G_1)+l^*(G_2)\le l^*(G)+1$.  \qed


We can now present a class of examples where the gap between the proper leafage
and the maximum number of MEP's in induced subgraphs becomes arbitrarily large.
The $n$-{\it kite} $G_n$ is the graph with $3n+1$ vertices consisting of a path
$P$ on successive vertices $\VEC v0n$, plus vertices $\{\VEC u1n\}$ and
$\{\VEC w1n\}$ such that $N(u_i)=N(w_i)=\{v_{i-1},v_i\}$.  Although no induced
subgraph has more than four MEP's, $l^*(G_n)=n+2$ for $n\ge 2$.  Note first that
$l^*(G_1)=2$.  For $n=2$, the vertices $\{u_1,w_1,v_1,u_2,w_2\}$ induce
$K_{1,4}$, a subgraph in which the four leaves are MEP's.
Hence $l^*(G_2)\ge 4$.  For $n>2$, we form $G_n$ by identifying $v_{n-1}$ from
$G_{n-1}$ with $v_0$ from $G_1$.  By the second part of Lemma 15, we thus have
$l^*(G_n)\ge l^*(G_{n-1})+1$, so $l^*(G_n)\ge n+2$ by induction.
Fig.~3 illustrates a construction that achieves equality.  Note that all $2n$
simplicial vertices are EP's, but we have leaves in the proper subtree \rep\ for
only half of them.

\gpic{
\expandafter\ifx\csname graph\endcsname\relax \csname newbox\endcsname\graph\fi
\expandafter\ifx\csname graphtemp\endcsname\relax \csname newdimen\endcsname\graphtemp\fi
\setbox\graph=\vtop{\vskip 0pt\hbox{%
    \special{pn 8}%
    \special{pa 76 610}%
    \special{pa 331 610}%
    \special{pa 331 356}%
    \special{fp}%
    \graphtemp=.5ex\advance\graphtemp by 0.508in
    \rlap{\kern 0.102in\lower\graphtemp\hbox to 0pt{\hss $v_0$\hss}}%
    \special{pa 356 356}%
    \special{pa 356 102}%
    \special{fp}%
    \special{pa 102 636}%
    \special{pa 610 636}%
    \special{fp}%
    \special{pa 381 356}%
    \special{pa 381 610}%
    \special{pa 1093 610}%
    \special{pa 1093 356}%
    \special{fp}%
    \graphtemp=.5ex\advance\graphtemp by 0.000in
    \rlap{\kern 0.356in\lower\graphtemp\hbox to 0pt{\hss $u_1$\hss}}%
    \graphtemp=.5ex\advance\graphtemp by 0.712in
    \rlap{\kern 0.356in\lower\graphtemp\hbox to 0pt{\hss $w_1$\hss}}%
    \graphtemp=.5ex\advance\graphtemp by 0.508in
    \rlap{\kern 0.737in\lower\graphtemp\hbox to 0pt{\hss $v_1$\hss}}%
    \special{pa 1119 356}%
    \special{pa 1119 102}%
    \special{fp}%
    \special{pa 864 636}%
    \special{pa 1373 636}%
    \special{fp}%
    \special{pa 1144 356}%
    \special{pa 1144 610}%
    \special{pa 1856 610}%
    \special{pa 1856 356}%
    \special{fp}%
    \graphtemp=.5ex\advance\graphtemp by 0.000in
    \rlap{\kern 1.119in\lower\graphtemp\hbox to 0pt{\hss $u_2$\hss}}%
    \graphtemp=.5ex\advance\graphtemp by 0.712in
    \rlap{\kern 1.119in\lower\graphtemp\hbox to 0pt{\hss $w_2$\hss}}%
    \graphtemp=.5ex\advance\graphtemp by 0.508in
    \rlap{\kern 1.500in\lower\graphtemp\hbox to 0pt{\hss $v_2$\hss}}%
    \special{pa 1881 356}%
    \special{pa 1881 102}%
    \special{fp}%
    \special{pa 1627 636}%
    \special{pa 2136 636}%
    \special{fp}%
    \special{pa 1907 356}%
    \special{pa 1907 610}%
    \special{pa 2619 610}%
    \special{pa 2619 356}%
    \special{fp}%
    \graphtemp=.5ex\advance\graphtemp by 0.000in
    \rlap{\kern 1.881in\lower\graphtemp\hbox to 0pt{\hss $u_3$\hss}}%
    \graphtemp=.5ex\advance\graphtemp by 0.712in
    \rlap{\kern 1.881in\lower\graphtemp\hbox to 0pt{\hss $w_3$\hss}}%
    \graphtemp=.5ex\advance\graphtemp by 0.508in
    \rlap{\kern 2.263in\lower\graphtemp\hbox to 0pt{\hss $v_3$\hss}}%
    \special{pa 2644 356}%
    \special{pa 2644 102}%
    \special{fp}%
    \special{pa 2390 636}%
    \special{pa 2898 636}%
    \special{fp}%
    \special{pa 2669 356}%
    \special{pa 2669 610}%
    \special{pa 2924 610}%
    \special{fp}%
    \graphtemp=.5ex\advance\graphtemp by 0.000in
    \rlap{\kern 2.644in\lower\graphtemp\hbox to 0pt{\hss $u_4$\hss}}%
    \graphtemp=.5ex\advance\graphtemp by 0.712in
    \rlap{\kern 2.644in\lower\graphtemp\hbox to 0pt{\hss $w_4$\hss}}%
    \graphtemp=.5ex\advance\graphtemp by 0.508in
    \rlap{\kern 2.898in\lower\graphtemp\hbox to 0pt{\hss $v_4$\hss}}%
    \hbox{\vrule depth0.712in width0pt height 0pt}%
    \kern 3.000in
  }%
}%
}
\ce{Fig.~3.  Optimal proper subtree \rep\ of $G_4$.}

\sp
Finally, we show that proper leafage equals leafage for $\Koe$-free chordal
graphs.
\skipit{Although Roberts proved that the concepts of MEP and
EP are the same for these graphs, that does not immediately imply this result,
since we have seen that the number of EP's need not be an upper bound for the
proper leafage of a chordal graph.}

\TH {11}.
If $G$ is a $\Koe$-free non-clique chordal graph, then $l(G)$, $l^*(G)$,
$a(G)$, and the number of inequivalent MEP's in $G$ are equal.
\PF
Recall that $a(G)$ denotes the maximum size of an asteroidal set in $G$, and
let $m(G)$ be the number of inequivalent MEP's in $G$.  If any of
$a(G),m(G),l(G)$ is 2, then $G$ is an interval graph.  Being $\Koe$-free, such
a graph $G$ is also a unit-interval graph, and hence all the parameters equal 2.
Hence we may assume the values exceed 2.  We prove that
$m(G)= a(G)\le l(G)\le l^*(G)\le a(G)$.
By Theorem 1 and the definition of proper leafage, we need only prove the
first equality and the last inequality.

In proving that $m(G)=a(G)$, we may assume that $G$ is reduced.  Let $X$ be the
set of MEP's in $G$.  To prove that $m(G)\le a(G)$, we prove that $X$ is an
asteroidal set.  It suffices to show that $G-N[x]$ is connected for each
$x\in X$.  If not, then by the definition of
MEP, $G-N[x]$ has two components $H_1,H_2$ such that $U_1=U_2$, where
$U_i=\{y\in N(x)\st N(y)\cap V(H_i)\ne\nul\}$.  Choose $w\in U_1$,
$y\in N(w)\cap V(H_1)$, and $z\in N(y)\cap V(H_2)$; now $w,x,y,z$ induce the
forbidden $\Koe$.

We next prove that $a(G)\le m(G)$.  Among all maximum-sized asteroidal sets, let
$X$ be one that maximizes the sum $\sg(X)$ of the pairwise distances between the
elements.  By the argument in Lemma 7, we may assume that $X\esub S(G)$.
Consider $x\in X$.  In order to prove that $x$ is an MEP, it suffices to prove
that $G-N[x]$ has only one component.  Since vertices of $X-x$ are linked by
paths avoiding $N[x]$, all of $X-x$ belongs to the same component $H_1$ of
$G-N[x]$.  If there is another component $H_2$ of $G-N[x]$, a common neighbor
of $H_1$ and $H_2$ in $N(x)$ would create an induced $\Koe$.  Hence we could
replace $x$ in $X$ by a vertex of $H_2$ to obtain a maximum asteroidal set with
greater distance sum.

For the last inequality, we build a proper subtree \rep\ of $G$ with $a(G)$
leaves, by induction on $n(G)$, in which the vertices of a maximum asteroidal
set $X$ maximizing $\sg(X)$ appear at the leaves of the host tree.  Again we may
assume that $G$ is reduced.  As proved above, each $x\in X$ is an MEP.  Since
$G$ is $\Koe$-free, distinct components of $G-N[x]$ cannot have common neighbors
in $N(x)$.  By Remark 1, the vertices of $N(x)$ form a minimal cutset.
Hence $G-N[x]$ is connected.  Thus $G-N[x]$ cannot have an asteroidal set of
size $a(G)$, since $x$ would augment it to a larger asteroidal set in $G$.

Now the induction hypothesis guarantees a proper subtree \rep\ of $G-x$ with
$a(G-x)$ leaves.  If $a(G-x)=a(G)$, then every maximum asteroidal set in $G-x$
has a vertex $y \in N(x)$, including sets that maximize $\sg$ and have all
elements simplicial.  In such a \rep, we have $y$ at a leaf $q$ of the host.
We add a leaf adjacent to $q$, and to $x$ we assign the minimal subtree
containing $q$ and at least one vertex of the subtree representing each vertex
of $N(x)$.  Since $x$ is simplicial, $N_G(x)\esub N_{G-x}[y]$.  Since the
\rep\ of $G-x$ is proper, the result is a proper subtree \rep\ of $G$.

Finally, if $a(G-x)<a(G)$, then we take the proper subtree \rep\ of $G-x$
guaranteed by induction for $X-x$, find a vertex of the host at which all 
of $N(x)$ appear, extend those subtrees from that vertex to a new leaf $q$,
and add another edge from $q$ to another new vertex $q'$, assigning
$f(x)=\{q,q'\}$.  \qed

\SH
{References}
\frenchspacing
\BP[1]
H. Broersma, T. Kloks, D. Kratsch, and H. M\"uller,
Independent sets in asteroidal triple-free graphs,
{\it Proceedings of ICALP'97},
P. Degano, R. Gorrieri, A. Marchetti-Spaccamela, (eds.),
(Springer-Verlag, 1997), {\it Lect. Notes Comp. Sci.} 1256, 760--770.
\BP[2]
H. Broersma, T. Kloks, D. Kratsch, and H. M\"uller,
A generalization of AT-free graphs and a generic algorithm
for solving triangulation problems,
Memorandum No.~1385, University of Twente, Enschede,
The Netherlands, 1997.
\BP[3]
P.A.~Buneman, A characterization of rigid circuit graphs, \DM\ 9(1974),
205--212.
\BP[4]
R.P.~Dilworth, A decomposition theorem for partially ordered sets,
{\it Ann.\ Math.\ }51(1950), 161--166.
\BP[5]
R.P.~Dilworth, Some combinatorial problems on partially ordered sets,
{\it Combinatorial Analysis} (Bellman and Hall, eds.) {\it Proc.\ Symp.\ Appl.\
Math.\ }(Amer.\ Math.\ Soc 1960), 85--90.
\BP[6]
G.A.~Dirac, On rigid circuit graphs, {\it Abh.\ Math.\ Sem.\ Univ.\ Hamburg}
25(1961), 71--76.
\BP[7]
D.R.~Fulkerson and O.A.~Gross, Incidence matrices and interval graphs,
{\it Pac. J. Math.\ }15(1965), 835--855.
\BP[8]
F.~Gavril, The intersection graphs of subtrees in trees are exactly the chordal
graphs, \JCTB\ 16(1974), 47--56.
\BP[9]
F.~Gavril, Generating the maximum spanning trees of a weighted graph.
{\it J. Algorithms} 8(1987), 592--597.
\BP[10]
P.C.~Gilmore and A.J.~Hoffman, A characterization of comparability graphs and
of interval graphs.  {\it Canad. J. Math.\ }16(1964), 539--548.
\BP[11]
M.C.~Golumbic, {\it Algorithmic Graph Theory and Perfect Graphs},
(Academic Press 1980).
\BP[12]
T. Kloks, D. Kratsch, and H. M\"uller, Asteroidal sets in graphs,
{\it Proceedings of WG'97}, R.~M\"ohring, (ed.), (Springer-Verlag, 1997),
{\it Lect. Notes Comp. Sci.} 1335, 229--241.
\BP[13]
T. Kloks, D. Kratsch, and H. M\"uller, On the structure of graphs with bounded
asteroidal number, Forschungsergebnisse Math/Inf/97/22, FSU Jena, Germany, 1997.
\BP[14]
B.~Leclerc, Arbres et dimension des ordres, \DM\ 14(1976), 69--76.
\BP[15]
C.B.~Lekkerkerker and J.Ch.~Boland, Representation of a finite graph by a set
of intervals on the real line.  {\it Fund.\ Math.\ }51(1962), 45--64.
\BP[16]
I.-J.~Lin, M.K.~Sen, and D.B.~West, Leafage of directed graphs, to appear.
\BP[17]
T.A.~McKee, Subtree catch graphs, {\it Congr.\ Numer.\ } 90(1992), 231--238.
\BP[18]
E.~Prisner, Representing triangulated graphs in stars,
{\it Abh.\ Math.\ Sem.\ Univ.\ Hamburg} 62(1992), 29--41.
\BP[19]
F.S.~Roberts, Indifference graphs, {\it Proof Techniques in Graph Theory}
(F.~Harary, ed.).  Academic Press (1969), 139--146.
\BP[20]
D.J.~Rose, Triangulated graphs and the elimination process,
{\it J. Math.\ Ann.\ Appl.\ }32(1970), 597--609.
\BP[21]
D.J.~Rose, R.E.~Tarjan, and G.S.~Leuker, Algorithmic aspects of vertex
elimination on graphs, {\it SIAM J. Comp.\ }5(1976), 266--283.
\BP[22]
Y.~Shibata, On the tree representation of chordal graphs. \JGT\ 12(1988), 
421--428.
\BP[23]
J.R.~Walter, Representations of chordal graphs as subtrees of a tree,
\JGT\  2(1978), 265--267.
\bye